\documentclass[reqno]{amsart}
\usepackage{amssymb}
\usepackage{graphicx}
\usepackage{pictex}
\usepackage{mathrsfs}
\usepackage{hyperref}

\begin{document}

\newtheorem{thm}{Theorem}
\newtheorem{lem}[thm]{Lemma}
\newtheorem{cor}[thm]{Corollary}
\newtheorem{prop}[thm]{Proposition}
\newtheorem{propa}[thm]{Proposition}

\theoremstyle{definition}
\newtheorem{defn}[thm]{{\color{green}Definition}}

\theoremstyle{remark}
\newtheorem{rmk}[thm]{Remark}

\newcommand{\CPtwo}{\mathbf{CP}^2}
\newcommand{\SLtwoC}{\mathrm{SL}(2,{\mathbf C})}
\newcommand{\SLtwoR}{\mathrm{SL}(2,{\mathbf R})}
\newcommand{\SLtwoZ}{\mathrm{SL}(2,{\mathbf Z})}
\newcommand{\QQ}{\hat{\mathbf Q} \cap [0,\infty]}

\newcommand{\nn}{\noindent}

\newenvironment{pf}{\noindent {\sc Proof.}\,\,}{\qed \vskip 6pt}
\newenvironment{pf*}{\noindent {\sc Proof}\,\,}{\qed \vskip 6pt}

\title[Uniqueness of Markoff Numbers which are Prime Powers]%
{An elementary proof of uniqueness of Markoff numbers which are prime powers}%
\author{Ying Zhang}
\address{Department of Mathematics, Yangzhou University, Yangzhou 225002, CHINA} %
\email{yingzhang@yzu.edu.cn}
\address{Current address: IMPA, Estrada Dona Castorina 110, 22460 Rio de Janeiro, BRAZIL}
\email{yiing@impa.br}
\date{Version 1: June 9, 2006. Version 2: January 31, 2007}
\thanks{Supported by a CNPq-TWAS postdoctoral
fellowship and in part by NSFC grant \#10671171.}
\begin{abstract}
We present a very elementary proof of the uniqueness of Markoff
numbers which are prime powers or twice prime powers, in the sense
that it uses neither algebraic number theory nor hyperbolic
geometry.
\end{abstract}

\maketitle

\section{Introduction}\label{s:intro}

\subsection{Markoff numbers} In his celebrated work on the minima of
indefinite binary quadratic forms, A. A. Markoff
\cite{markoff1880ma} was naturally led to the study of Diophantine
equation---now known as the {\it Markoff equation}
\begin{eqnarray}\label{eqn:markoff}
x^2+y^2+z^2=3xyz.
\end{eqnarray}
The solution triples $(x,y,z)$ in positive integers are called by
Frobenius \cite{frobenius1913} the {\it Markoff triples}, and the
individual positive integers occur the {\it Markoff numbers}.

For convenience, we shall not distinguish a Markoff triple from its
permutation class, and when convenient, usually arrange its elements
in ascending order. Following Cassels \cite{cassels1957book}, we
call the Markoff triples $(1,1,1)$ and $(1,1,2)$ {\it singular},
while all the others {\it non-singular}. It is easy to show that the
elements of a non-singular Markoff triple are all distinct.

In ascending order of their maximal elements, the first 12 Markoff
triples are:
\begin{align*}
&(1,1,1),\,(1,1,2),\,(1,2,5),\,(1,5,13),\,(2,5,29),\,(1,13,34),\,(1,34,89), \\
&(2,29,169),\,(5,13,194),\,(1, 89, 233),\,(5, 29, 433),\,(89, 233, 610); %
\end{align*}
while the first 40 Markoff numbers as recorded in \cite{sloane} are:
\begin{align*}
& 1,\,2,\,5,\,13,\,29,\,34,\,89,\,169,\,194,\,233,\,433, %
610, 985, 1325, 1597, 2897, \\ & 4181, 5741, 6466, 7561, 9077, 10946, 14701, %
28657, 33461, 37666, 43261, \\ & 51641, 62210, 75025, 96557, 135137, %
195025, 196418, 294685, 426389, \\ & 499393, 514229, 646018, 925765. %
\end{align*}

\subsection{Sketch of Markoff's work}\label{ss:mwork}

Let $f(\xi,\eta)=a\xi^2+b\xi\eta+c\eta^2$ be a binary quadratic form
with real coefficients. The discriminant of $f$ is defined as
$\delta(f)=b^2-4ac$, and the minimum $m(f)$ of $f$ is defined as
$$m(f)=\inf |f(\xi,\eta)|,$$
where the infimum is taken over all pairs of integers $\xi,\eta$ not
both zero.

Two quadratic forms $f(\xi,\eta)$ and $g(\xi,\eta)$ are said to be
{\it equivalent} if there exist integers $a,b,c,d$ such that
$ad-bc=\pm 1$ and $f(a\xi+b\eta,c\xi+d\eta)=g(\xi,\eta)$.

Then Markoff's aforementioned work on the minima of real indefinite
binary quadratic forms can be stated as follows.

\vskip 6pt

{\sc Markoff's Theorem.}\, {\it Let $f$ be a real indefinite binary
quadratic form. Then inequality $m(f)/\sqrt{\delta(f)}>1/3$ holds if
and only if $f$ is equivalent to a multiple of a Markoff form.}

\vskip 6pt

Here the {\it Markoff form} associated to a Markoff triple
$(m,m_1,m_2)$ with $m \ge m_1 \ge m_2$ is defined as an indefinite
binary quadratic form with integer coefficients, as follows.
First, let $u$ be the least non-negative integer such that
\begin{eqnarray*}
um_1\equiv m_2 \quad ({\rm mod}\; m) \quad \text{or} \quad %
um_1\equiv -m_2 \quad ({\rm mod}\; m).
\end{eqnarray*}
Since $0 \equiv m(3m_1m_2-m)=m_1^2+m_2^2\equiv (u^2+1)m_1^2 \;({\rm
mod}\; m)$ and, as will be shown in \S \ref{ss:basic},
$\gcd(m_1,m)=1$, we have $u^2+1 \equiv 0 \:({\rm mod}\; m)$. Now let
\begin{eqnarray*}v=(u^2+1)/m.\end{eqnarray*} %
The Markoff form associated to Markoff triple $(m,m_1,m_2)$ is then
defined as
\begin{eqnarray}
\phi_{(m,m_1,m_2)}(\xi,\eta)=m\xi^2+(3m-2u)\xi\eta+(v-3u)\eta^2.
\end{eqnarray}
Note that for $\phi:=\phi_{(m,m_1,m_2)}$ we have
$\delta(\phi)=9m^2-4$ and $m(\phi)=m$. %

\vskip 6pt

\nn {\bf Remark.}\, Note that Markoff \cite{markoff1879ma}
\cite{markoff1880ma} used continued fractions to obtain his results,
and his proofs were only sketched. Dickson \cite[Ch.VII]{dickson}
gave a detailed interpretation of it. Frobenius \cite{frobenius1913}
made a systematic study of the Markoff numbers, based on which Remak
\cite{remak1924ma} presented a proof of Markoff's Theorem using no
continued fractions. Markoff's above result also has a well-known
equivalent formulation in terms of the approximation of irrationals
by rationals; see Cassels \cite{cassels1957book} and Cusick--Flahive
\cite{cusick-flahive} for detailed explanations.

\subsection{Neighbors of a Markoff triple}

\begin{figure}
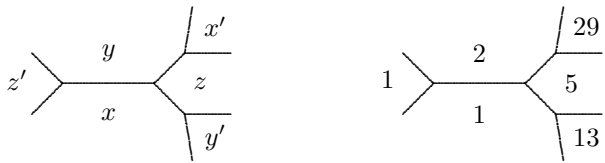

\begin{center}
\mbox{\beginpicture \setcoordinatesystem units <0.04in,0.04in>
\setplotarea x from -20 to 20 , y from -10 to 10
\plot -6 0 6 0 / %
\plot -10 4 -6 0 -10 -4 / \plot 10 4 6 0 10 -4 / %
\plot 11 10 10 4 16 4 / \plot 11 -10 10 -4 16 -4 / %
\put {\mbox{\normalsize $y$}}  [cc] <0mm,0mm>  at  0  4 %
\put {\mbox{\normalsize $x$}}  [cc] <0mm,0mm>  at  0 -4 %
\put {\mbox{\normalsize $z$}}  [cc] <0mm,0mm>  at  12  0 %
\put {\mbox{\normalsize $y'$}} [cc] <0mm,0mm>  at  14 -7 %
\put {\mbox{\normalsize $x'$}} [cc] <0mm,0mm>  at  14 7.5 %
\put {\mbox{\normalsize $z'$}} [cc] <0mm,0mm>  at -12 0.5 %
\endpicture}
\hspace{0.25in} \mbox{\beginpicture \setcoordinatesystem units
<0.04in,0.04in> \setplotarea x from -20 to 20 , y from -10 to 10
\plot -6 0 6 0 / %
\plot -10 4 -6 0 -10 -4 / \plot 10 4 6 0 10 -4 / %
\plot 11 10 10 4 16 4 / \plot 11 -10 10 -4 16 -4 / %
\put {\mbox{\normalsize $2$}}  [cc] <0mm,0mm>  at  0  4 %
\put {\mbox{\normalsize $1$}}  [cc] <0mm,0mm>  at  0 -4 %
\put {\mbox{\normalsize $5$}}  [cc] <0mm,0mm>  at  12  0 %
\put {\mbox{\normalsize $13$}} [cc] <0mm,0mm>  at  14 -7 %
\put {\mbox{\normalsize $29$}} [cc] <0mm,0mm>  at  14 7.5 %
\put {\mbox{\normalsize $1$}} [cc] <0mm,0mm>  at -12 0.5 %
\endpicture}
\end{center}
\caption{Tree structure of a Markoff triple and its neighbors}\label{fig:3neighbors} %
\end{figure}

That Markoff equation (\ref{eqn:markoff}) is particularly
interesting lies in the fact that it is a quadratic equation in each
of the variables, and hence new solutions can be obtained by a
simple process from a given one, $(x,y,z)$.
To see this, keep $x$ and $y$ fixed and let $z'$ be the other root
of (\ref{eqn:markoff}), regarded as a quadratic equation in $z$.
Rewriting (\ref{eqn:markoff}) as $z^2-3xyz+(x^2+y^2)=0$, we have
$z+z'=3xy$ and $zz'=x^2+y^2$. Thus $z'$ is a positive integer and
$(x,y,z')$ is another solution triple to (\ref{eqn:markoff}) in
positive integers, that is, a Markoff triple. Similarly, we obtain
two other Markoff triples $(x',y,z)$ and $(x,y',z)$. We call these
three new Markoff triples thus obtained the {\it neighbors} of the
$(x,y,z)$. See Figure \ref{fig:3neighbors} for an illustration.

\subsection{Reduction}

In \cite[pp.397--398]{markoff1880ma}, Markoff showed that every
Markoff triple can be obtained from the simplest by appropriately
iterating the above process.

\vskip 6pt

{\sc The Reduction Theorem.}  {\it Every Markoff triple can be
traced back to $(1,1,1)$ by repeatedly performing the following
operation on Markoff triples:
\begin{eqnarray}\label{eqn:operation}
(x,y,z)\longmapsto(x,y,z'):=(x,y,3xy-z), %
\end{eqnarray}
where the elements of $(x,y,z)$ is arranged so that $x \le y \le
z$.}

\vskip 6pt

Note that to perform the next operation, one needs to first
rearrange the elements of $(x,y,z')$ in ascending order. As an
example, we see
\begin{eqnarray*}
& (13,194,7561) \longmapsto (13,194,5)\sim(5,13,194) \longmapsto (5,13,1)\sim(1,5,13)  \\ %
& \longmapsto (1,5,2)\sim(1,2,5) \longmapsto (1,2,1)\sim(1,1,2) \longmapsto (1,1,1). %
\end{eqnarray*}

A simple proof of the theorem is given in
\cite[pp.27--28]{cassels1957book}; see also
\cite[pp.17--18]{cusick-flahive}. The idea is that operation
(\ref{eqn:operation}) reduces the maximal elements of Markoff
triples as long as the input triple is non-singular. Indeed, one has
$x<y<z$ and $(z-y)(z'-y)=zz'-(z+z')y+y^2=x^2+2y^2-3xy^2<0$; hence
$z'<y$.

Here we give a slightly different proof, the idea of which we get
from \cite{bowditch1998plms}.

\vskip 6pt

\noindent {\sc Proof.} \, The operation $(x,y,z)\longmapsto(x,y,z')$
reduces the lengths, $x+y+z$, of Markoff triples exactly when
$z'<z$. Therefore, after a finite number of times of length
reduction, one stops when $z' \ge z$, or equivalently, $2z \le 3xy$.
We claim that $z=1$ in this case, and hence $(x,y,z)=(1,1,1)$.
Indeed, if $z \ge 2$ then one obtains from $x \le y \le z$ and $2z
\le 3xy$ that
\begin{eqnarray}\label{eqn:632}
1=\frac{x}{3yz}+\frac{y}{3zx}+\frac{z}{3xy} \le %
\frac{1}{6}+\frac{1}{3}+\frac{1}{2}=1.
\end{eqnarray}
This forces that $x=y,z=2$ and $x=1,y=z$ both hold, a contradiction.
\qed

\vskip 6pt

\begin{figure}
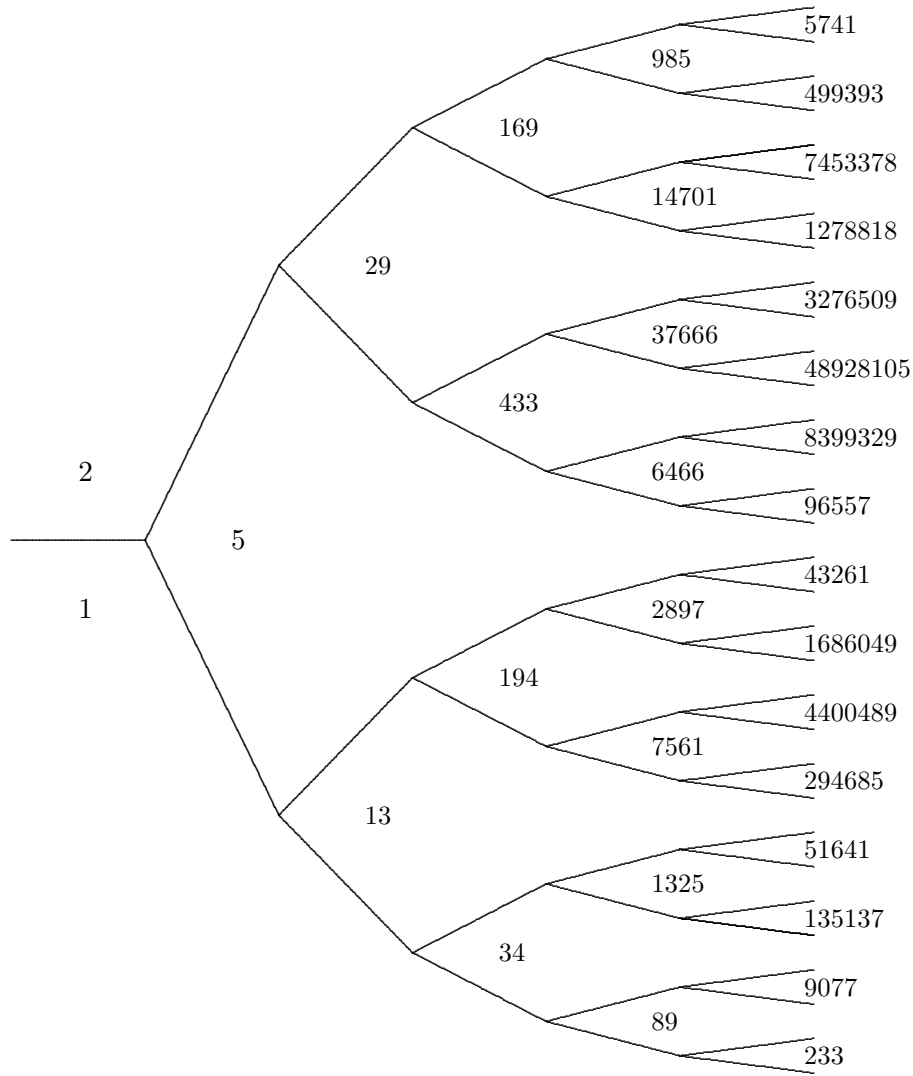

\begin{center}
\mbox{\beginpicture \setcoordinatesystem units <0.01in,0.018in>
\setplotarea x from 0 to 480 , y from -160 to 160
\plot 0 0 70 0 140 80 210 120 280 140 350 150 420 155 / %
\plot 70 -0 140 -80 210 -120 280 -140 350 -150 420 -155 / %
\plot 140 80 210 40 280 20 350 10 420 5 / %
\plot 140 -80 210 -40 280 -20 350 -10 420 -5 / %
\plot 210 120 280 100 350 90 420 85 / %
\plot 210 -120 280 -100 350 -90 420 -85 / %
\plot 210 40 280 60 350 70 420 75 / %
\plot 210 -40 280 -60 350 -70 420 -75 / %
\plot 280 140 350 130 420 125 / %
\plot 280 -140 350 -130 420 -125 / %
\plot 280 100 350 110 420 115 / %
\plot 280 -100 350 -110 420 -115 / %
\plot 280 60 350 50 420 45 / %
\plot 280 -60 350 -50 420 -45 / %
\plot 280 20 350 30 420 35 / %
\plot 280 -20 350 -30 420 -35 / %
\plot 350 150 420 145 / \plot 350 -150 420 -145 / %
\plot 350 130 420 135 / \plot 350 -130 420 -135 / %
\plot 350 110 420 115 / \plot 350 -110 420 -115 / %
\plot 350 110 420 105 / \plot 350 -110 420 -105 / %
\plot 350  90 420  95 / \plot 350  -90 420  -95 / %
\plot 350  70 420  65 / \plot 350  -70 420  -65 / %
\plot 350  50 420  55 / \plot 350  -50 420  -55 / %
\plot 350  30 420  25 / \plot 350  -30 420  -25 / %
\plot 350  10 420  15 / \plot 350  -10 420  -15 / %
\put {\mbox{\large $1$}}  [lc] <0mm,0mm>  at  35 -20 %
\put {\mbox{\large $2$}}  [lc] <0mm,0mm>  at  35 20 %
\put {\mbox{\large $5$}}  [lc] <0mm,0mm>  at 115 0 %
\put {\mbox{\normalsize $29$}}  [lc] <0mm,0mm>  at 185 80 %
\put {\mbox{\normalsize $13$}}  [lc] <0mm,0mm>  at 185 -80 %
\put {\mbox{\normalsize $169$}} [lc] <0mm,0mm>  at 255 120 %
\put {\mbox{\normalsize $433$}} [lc] <0mm,0mm>  at 255 40 %
\put {\mbox{\normalsize $194$}} [lc] <0mm,0mm>  at 255 -40 %
\put {\mbox{\normalsize $34$}} [lc] <0mm,0mm>  at 255 -120 %
\put {\mbox{\normalsize $985$}}  [lc] <0mm,0mm>  at 335 140 %
\put {\mbox{\normalsize $14701$}}  [lc] <0mm,0mm>  at 335 100 %
\put {\mbox{\normalsize $37666$}}  [lc] <0mm,0mm>  at 335 60 %
\put {\mbox{\normalsize $6466$}}  [lc] <0mm,0mm>  at 335 20 %
\put {\mbox{\normalsize $2897$}}  [lc] <0mm,0mm>  at 335 -20 %
\put {\mbox{\normalsize $7561$}}  [lc] <0mm,0mm>  at 335 -60 %
\put {\mbox{\normalsize $1325$}}  [lc] <0mm,0mm>  at 335 -100 %
\put {\mbox{\normalsize $89$}}  [lc] <0mm,0mm>  at 335 -140 %
\put {\mbox{\normalsize $5741$}}  [lc] <0mm,0mm>  at 415 150 %
\put {\mbox{\normalsize $499393$}}  [lc] <0mm,0mm>  at 415 130 %
\put {\mbox{\normalsize $7453378$}}  [lc] <0mm,0mm>  at 415 110 %
\put {\mbox{\normalsize $1278818$}}  [lc] <0mm,0mm>  at 415 90 %
\put {\mbox{\normalsize $3276509$}}  [lc] <0mm,0mm>  at 415 70 %
\put {\mbox{\normalsize $48928105$}}  [lc] <0mm,0mm>  at 415 50 %
\put {\mbox{\normalsize $8399329$}}  [lc] <0mm,0mm>  at 415 30 %
\put {\mbox{\normalsize $96557$}}  [lc] <0mm,0mm>  at 415 10 %
\put {\mbox{\normalsize $43261$}}  [lc] <0mm,0mm>  at 415 -10 %
\put {\mbox{\normalsize $1686049$}}  [lc] <0mm,0mm>  at 415 -30 %
\put {\mbox{\normalsize $4400489$}}  [lc] <0mm,0mm>  at 415 -50 %
\put {\mbox{\normalsize $294685$}}  [lc] <0mm,0mm>  at 415 -70 %
\put {\mbox{\normalsize $51641$}}  [lc] <0mm,0mm>  at 415 -90 %
\put {\mbox{\normalsize $135137$}}  [lc] <0mm,0mm>  at 415 -110 %
\put {\mbox{\normalsize $9077$}}  [lc] <0mm,0mm>  at 415 -130 %
\put {\mbox{\normalsize $233$}}  [lc] <0mm,0mm>  at 415 -150 %
\endpicture}
\end{center}
\caption{Markoff numbers in an infinite binary tree}\label{fig:xyz-tree} %
\end{figure}

\subsection{First properties of Markoff numbers}\label{ss:basic}

As an immediate corollary of the Reduction Theorem, we see that the
elements of a Markoff triple are pairwise coprime. Moreover, since
$zz'=x^2+y^2$ and $\gcd (x,y)=1$, a Markoff number is not a multiple
of $4$, and each odd prime factor of a Markoff number is $\equiv 1
\,({\rm mod}\; 4)$. Consequently, every odd Markoff number is
$\equiv 1 \,({\rm mod}\; 4)$ and every even Markoff number is
$\equiv 2 \,({\rm mod}\,\, 8)$. Indeed, it is shown in \cite{zy}
that every even Markoff number is $\equiv 2 \,({\rm mod}\,\, 32)$.

\subsection{An illustration}

The Reduction Theorem tells that, starting from $(1,2,5)$ and
generating new neighbors repeatedly, one will obtain all the Markoff
triples. This is depicted as an infinite binary tree in Figure
\ref{fig:xyz-tree} in which all the Markoff numbers appear in the
regions while all non-singular Markoff triples appear around
vertices. In this shape it seems to be first drawn by Thomas E. Ace
on his web-page \verb|http://www.minortriad.com/markoff.html|.

\subsection{The uniqueness problem}

A problem then arises naturally: \emph{Does every Markoff number
appear exactly once in the regions in Figure \ref{fig:xyz-tree}?} In
other words, are there any repetitions among all the numbers occur?

The following conjecture on the uniqueness of Markoff
numbers/triples was first mentioned explicitly by G. Frobenius as a
question in his 1913 paper \cite{frobenius1913}. It asserts that a
Markoff triple is uniquely determined by its maximal element. (And
we shall simply say that a Markoff number $z$ is {\it unique} if the
following is true for $z$.)

\vskip 6pt

{\sc The Unicity Conjecture.} \, {\it Suppose \,$(x,y,z)$ and
$(\tilde{x},\tilde{y},z)$ are Markoff triples with $x \le y \le z$
and \,$\tilde{x} \le \tilde{y} \le z$. Then $x=\tilde{x}$ and
$y=\tilde{y}$.}

\vskip 6pt

The conjecture has been proved only for some special subsets of the
Markoff numbers. The following affirmative result for Markoff
numbers which are prime powers  or twice prime powers was first
proved independently and partly by Baragar \cite{baragar1996cmb},
Button \cite{button1998jlms} and Schmutz \cite{schmutz1996ma} using
either algebraic number theory
(\cite{baragar1996cmb},\cite{button1998jlms}) or hyperbolic geometry
(\cite{schmutz1996ma}). And a stronger result along the same lines
has been obtained later by Button in \cite{button2001jnt}; in
particular, a Markoff number is shown to be unique if it is a
``small'' ($\le 10^{35}$) multiple of a prime power.

\begin{thm}[Baragar \cite{baragar1996cmb}; Button
\cite{button1998jlms}; Schmutz \cite{schmutz1996ma}]\label{thm:button-schmutz} %
A Markoff number is unique if it is either a prime power or twice a
prime power.
\end{thm}

This paper is motivated by a simple proof of Theorem
\ref{thm:button-schmutz} recently published by Lang and Tan
\cite{lang-tan2005markoff}, which uses some elementary facts from
the hyperbolic geometry of the modular torus with one cusp, as used
by Cohn in \cite{cohn1955am}.
The aim of this paper is to present in detail a completely
elementary proof of Theorem \ref{thm:button-schmutz} that uses
neither algebraic number theory nor hyperbolic geometry so that an
average reader will be able to fully understand it with no
difficulty. Though it is later clear that all the needed ingredients
of the proof were already known as early as 1913 in Frobenius' work,
we must admit that we first obtain them from hyperbolic geometry,
especially, that used in \cite{lang-tan2005markoff} and
\cite{cohn1955am}.

The rest of the paper is organized as follows. In \S \ref{s:slopes}
we parametrize Markoff numbers using non-negative rationals (slopes)
$t\in\QQ$. We also define $u_t$ as in \S \ref{ss:mwork} and verify
some properties of the pairs $(m_t, u_t)$. Then in \S \ref{s:pf},
with the help of a simple lemma (Lemma \ref{lem:roots}), we give the
promised elementary proof of Theorem \ref{thm:button-schmutz}. In \S
\ref{s:M-matrices} we introduce the so-called Markoff matrices to
generate all Markoff numbers. Certain properties of these matrices
are then discussed in \S \ref{s:properties}. In particular,
alternative proofs of Lemmas \ref{lem:increasing} and \ref{lem:um}
will be given. Finally, in \S \ref{s:rmk} we give a geometric
explanation of the Markoff numbers and related numbers.

\vskip 6pt

\noindent {\bf Acknowledgements.} The author would like to thank Ser
Peow Tan for helpful conversations and suggestions. Thanks are also
due to a referee of the first version of this note, whose
constructive suggestions helped improve the exposition of the
current version.

\section{Slopes of Markoff numbers}\label{s:slopes}

\subsection{Slopes of Markoff numbers}

It is natural and very useful to associate to each Markoff number a
slope, that is, an ordered pair of non-negative coprime integers.
This was first done by Frobenius in \cite{frobenius1913} where he
set
\begin{eqnarray*}
m(1,0)=1,\;m(0,1)=2,\;m(1,1)=5,\;m(1,2)=13,\;m(2,1)=29,\; \ldots %
\end{eqnarray*}
(These pairs are also called by Cusick and Flahive
\cite{cusick-flahive} the {\sc Frobenius coordinates} of Markoff
numbers). Note that, by identifying $(\mu,\nu)$ with $\nu/\mu$, the
slopes are nothing but the positive rationals together with $0$ and
$\infty$. In the latter context we shall write\, $m_r$ for
$m(\mu,\nu)$ where $r=\nu/\mu$.

Let us write $\hat{\mathbf Q}:={\mathbf Q} \cup \{\infty\}$. We
shall also call $\infty=\frac10=\frac{-1}{0}$ a rational. Then the
set of slopes we consider is the set of rationals in $[0,\infty]$,
that is, $\QQ$.

\subsection{Farey sum of rationals}\label{ss:farey} %

\nn There is a simple but useful way to obtain all the positive
rationals by making the so-called Farey sums repeatedly.
Specifically, starting with $0=\frac01$ and $\infty=\frac10$ (of
level $0$), all positive rationals can be generated, level by level,
as follows:
\begin{eqnarray*}
&\frac11=\frac{0+1}{1+0}; \\ %
&\frac12=\frac{0+1}{1+1}, \frac21=\frac{1+1}{1+0}; \\ %
&\frac13=\frac{0+1}{1+2}, \frac23=\frac{1+1}{2+1},
 \frac32=\frac{1+2}{1+1}, \frac31=\frac{2+1}{1+0}; \\ %
&\frac14=\frac{0+1}{1+3}, \frac25=\frac{1+1}{3+2},
\frac35=\frac{1+2}{2+3}, \frac34=\frac{2+1}{3+1},
\frac43=\frac{1+3}{1+2}, \frac53=\frac{3+2}{2+1},
\frac52=\frac{2+3}{1+1}, \frac41=\frac{3+1}{1+0}; %
\end{eqnarray*}
and so on. (To obtain the negative rationals, one starts from
$\infty=\frac{-1}{0}$ and $0=\frac01$ instead and makes the Farey
sums recursively as above). In particular, we have the notion of
{\it Farey level} for positive rationals, with levels 1 to 4 shown
as above.
To obtain all the rationals in $[0,\infty]$ of Farey level $n+1$, we
simply start with all those of Farey level not exceeding $n$,
arrange them in ascending order, and make Farey sum for each pair of
adjacent ones among them. In particular, we are allowed to prove a
proposition concerning all the positive rationals by induction on
the Farey levels of the rationals involved. In what follows we shall
make the above idea precise and present some basic facts that will
be needed in later part of this paper.

By the standard reduced form of a rational number $t$ we mean the
unique fractional expression $t={\nu}/{\mu}$ where $\mu,\nu$ are
coprime integers with $\mu \ge 0$. Two rationals $r, s$ are said to
be {\it Farey neighbors} (and that they form a {\it Farey pair}) if
they have standard reduced forms $r=b/a$ and $s=d/c$ so that
$ad-bc=\pm 1$. Given a Farey pair $r, s$ with standard reduced forms
$r=b/a$ and $s=d/c$, their {\it Farey sum} is defined as
\begin{eqnarray}\label{eqn:fareyadd}
r\oplus s:=\frac{b+d}{a+c}
\end{eqnarray}
which is certainly in its standard reduced form. (Note that in terms
of $(a,b)$ and $(c,d)$ regarded as plane vectors, the Farey sum is
just the vector sum). Clearly, $r \oplus s = s \oplus r$. It is easy
to see that $r\oplus s$ falls in between $r$ and $s$ and is a common
Farey neighbor of $r$ and $s$. We shall call the ordered triple
$(r,t,s)$ a {\it Farey triple}.

It follows from the Euclidean algorithm that every positive rational
can be written in a unique way as the Farey sum of a Farey pair of
rationals in $[0,\infty]$. Indeed, for a given positive rational
$t$, among all its Farey neighbors there are exactly two, $r$ and
$s$, having smaller or the same denominators; it can be easily shown
that $r$ and $s$ form a Farey pair and $t = r \oplus s$. We call $r$
and $s$ the {\it direct descents} of $t$. As a consequence, it is
easy to see that in every Farey pair in $\hat{\mathbf Q} \cap
[0,\infty]$, the one with smaller denominator or numerator has
smaller Farey level and is a direct descent of the other. Hence it
can be shown by induction that all rationals between $\frac01$ and
$\frac10$ will appear in the above process of recursively making
Farey sums.

To end this subsection, we give a formal definition of the notion of
Farey level. First, we set the Farey level of each of $\frac01$ and
$\frac10$ to be $0$. Recursively, for a Farey pair $r,s$ in
$\hat{\mathbf Q} \cap [0,\infty]$, we define the Farey level of
their Farey sum $t=r \oplus s$ to be the sum of their respective
Farey levels. In this way we then have recursively defined a Farey
level for each $t \in \hat{\mathbf Q} \cap [0,\infty]$.

\begin{figure}
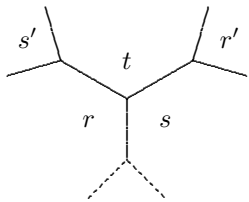

\begin{center}
\mbox{\beginpicture \setcoordinatesystem units <0.04in,0.04in>
\setplotarea x from -20 to 20 , y from -12 to 15
\plot 0 0 0 -8 / \plot -8.66 5 0 0 8.66 5 / %
\plot 15.66 3 8.66 5 10.66 12 / %
\plot -15.66 3 -8.66 5 -10.66 12 / %
\setdashes<1.5pt> %
\plot -5 -13 0 -8 5 -13 / %
\put {\mbox{\normalsize $t$}}  [cc] <0mm,0mm>  at  0  5 %
\put {\mbox{\normalsize $r$}}  [cc] <0mm,0mm>  at -5 -3 %
\put {\mbox{\normalsize $s$}}  [cc] <0mm,0mm>  at  5 -3 %
\put {\mbox{\normalsize $r'$}} [cc] <0mm,0mm>  at  13.5 8 %
\put {\mbox{\normalsize $s'$}} [cc] <0mm,0mm>  at -13 8 %
\endpicture}
\end{center}
\caption{Farey triples $(r,t,s)$, $(r,s',t)$ and $(t,r',s)$}\label{fig:triples} %
\end{figure}


\subsection{Further properties of Markoff numbers}\label{ss:furtherp} %


It is easy to see that, in terms of slopes, each Markoff triple is
then of the (more natural) form $(m_r,m_t,m_s)$ where $(r,t,s)$ is a
Farey triple in $\QQ$ with $r<t<s$.

We define for each $t\in \QQ$ an integer $u_t$ with $0 \le u_t \le
m_t$ as follows. First, we set $u_{0/1}=0$, $u_{1/0}=1$. In general,
for $t \in {\mathbf Q}\cap (0,\infty)$, $u_t$ is defined by
\begin{eqnarray}\label{eqn:u}
u_t \equiv m_s / m_r \; ({\rm mod}\; m_t).
\end{eqnarray}
Then $u_t$ depends only on $t$ but not the triple $(r,t,s)$ since
for the neighboring Farey triples $(r,s',t)$ and $(t,r',s)$ as shown
in Figure \ref{fig:triples} we have
\begin{eqnarray}\label{eqn:ratios}
m_r / m_{s'} \equiv m_s / m_r \equiv m_{r'} / m_s \; ({\rm mod}\; m_t) %
\end{eqnarray}
which in turn follows from
\begin{eqnarray*}
m_sm_{s'}=m_r^2+m_t^2 \quad\,\, \text{and} \quad\,\, m_rm_{r'}=m_t^2+m_s^2.
\end{eqnarray*}
Now since $0 \equiv m_r^2+m_s^2 \equiv m_r^2(1+u_t^2) \; ({\rm
mod}\; m_t)$ and $\gcd(m_r,m_t)=1$, we have from (\ref{eqn:u})
\begin{eqnarray}\label{eqn:uu}
u_t^2 +1 \equiv 0 \quad ({\rm mod}\; m_t). %
\end{eqnarray}

Furthermore, we have
\begin{lem}\label{lem:increasing}
The ratio $u_t/m_t$ is strictly increasing with respect to $t \in
\QQ$. In particular, $0 \le u_t \le m_t/2$, with strict inequalities
for $t \neq \frac01, \frac10$.
\end{lem}

In fact, Lemma \ref{lem:increasing} follows from the following %

\begin{lem}\label{lem:um} For every Farey triple $(r,t,s)$ in $\QQ$
with $r<t<s$,
\begin{eqnarray}\label{eqn:diff}
 \frac{u_t}{m_t}-\frac{u_r}{m_r}=\frac{m_s}{m_r m_t} %
 \quad\quad {\rm and} \quad\quad %
 \frac{u_s}{m_s}-\frac{u_t}{m_t}=\frac{m_r}{m_t m_s},
\end{eqnarray}
which are equivalent respectively to
\begin{eqnarray}\label{eqn:um}
 u_t m_r - u_r m_t = m_s \quad\quad {\rm and} \quad\quad %
 u_s m_t - u_t m_s = m_r.
\end{eqnarray}
\end{lem}

\begin{pf} We prove it by induction on the Farey levels of the
rationals involved. The conclusion is easily checked to be true for
the Farey triple $(\frac01,\frac11,\frac10)$. Now suppose that
(\ref{eqn:um}) holds for all Farey triples $(r,t,s)$ in $\QQ$ with
$r<t<s$ and with the Farey level of $t$ not exceeding $n \ge 1$. In
particular, this implies that $0 \le {u_r}/{m_r} < {u_t}/{m_t} <
{u_s}/{m_s} \le 1/2$.

Then we only need to show that (\ref{eqn:um}) also holds for the
Farey triples $(r,s',t)$ and $(t,r',s)$ as shown in Figure
\ref{fig:triples}. Since the proofs for the two cases are entirely
similar, we prove it for the case $(r,s',t)$ only, that is, we show
that
\begin{eqnarray}\label{eqn:ums'}
 u_{s'} m_r - u_r m_{s'} = m_t \quad\quad {\rm and} \quad\quad %
 u_t m_{s'} - u_{s'} m_t = m_r.
\end{eqnarray}
For this, we first see from (\ref{eqn:um}) and
$m_r^2+m_t^2=m_{s'}m_s$ that
\begin{eqnarray}\label{eqn:us'}
 u:=\frac{m_t + u_r m_{s'}}{m_r} =\frac{u_t m_{s'}-m_r}{m_t}.
\end{eqnarray}
Note that $0<u/m_{s'}<{u_t}/{m_t} < 1/2$ and, by (\ref{eqn:ratios}),
$u$ is an integer. Hence (\ref{eqn:ums'}) holds with $u_{s'}$
replaced by $u$. But this in turn implies that $u \equiv m_t/m_r
\;({\rm mod}\; m_{s'})$, and hence $u=u_{s'}$ by the definition of
$u_{s'}$. This proves Lemma \ref{lem:um}.
\end{pf}

\vskip 6pt

\nn {\bf Remark.}\, The inequalities in (\ref{eqn:um}) first
appeared in \cite[p.602]{frobenius1913}, though they were contained
essentially but implicitly in \cite{markoff1880ma}. The result of
Lemma \ref{lem:increasing} was stated and proved by Remak in
\cite{remak1924ma}.
In later part of this paper (see \S \ref{ss:pfum}), Lemma
\ref{lem:um} will also be obtained in a nice way as a corollary of
the properties (see Proposition \ref{prop:Mr-properties}) of the
so-called Markoff matrices which are interesting in their own right.

\subsection{Slope form of the Unicity Conjecture}

In terms of slopes, we may rephrase the Unicity Conjecture as

\vskip 6pt

{\sc The Unicity Conjecture (Slope form)}.\, {\it The Markoff
numbers $m_t$, $t \in \QQ$ are all distinct.}


\section{Proof of Theorem \ref{thm:button-schmutz}}\label{s:pf}


\nn We are now ready to give a very elementary proof for Theorem
\ref{thm:button-schmutz}, using Lemma \ref{lem:increasing} and the
following simple lemma whose proof can be found in \cite{zy}.

\begin{lem}\label{lem:roots}
Suppose $m=p^n$ or $2p^n$ for an odd prime $p$ and an integer $n \ge
1$. Then, for any integer $l$ coprime to $m$, the binomial
congruence equation
\begin{eqnarray}\label{eqn:qbinom}
x^2 + l \equiv 0 \quad ({\rm mod}\; m)
\end{eqnarray}
has at most one integer solution $x$ with $0<x<m/2$.
\end{lem}

\nn {\sc Proof of Theorem \ref{thm:button-schmutz}.}\; Suppose there
exist slopes $t,t^\ast\in \QQ$ such that
$$m_t=m_{t^\ast}=p^n \;\; \text{or} \;\; 2p^n$$ %
for an odd prime $p$ and an integer $n \ge 1$. By (\ref{eqn:uu}) and
its analog for $u_{t^\ast}$, Lemma \ref{lem:roots} applies to give
$u_t=u_{t^\ast}$. Then $t=t^\ast$ by Lemma \ref{lem:increasing}.
This proves Theorem \ref{thm:button-schmutz}. \qed

\vskip 6pt

\nn {\bf Remark.}\, The reader who is interested in merely the proof
of Theorem \ref{thm:button-schmutz} may well exit here. The rest of
this paper is devoted to a discussion of the so-called Markoff
matrices, which (with the exception of \S \ref{ss:pfum}) constitutes
the main body of an earlier version of this paper and can be used to
prove our earlier results in a nice way.


\section{Markoff matrices}\label{s:M-matrices}


\nn It is Harvey Cohn \cite{cohn1955am} who first noticed the
relationship of Markoff equation (\ref{eqn:markoff}) and one of
Fricke's trace identities, (\ref{eqn:-2}) below, for matrices in
$\SLtwoC$. This gives us a nice way to generate the Markoff numbers
using the so-called Markoff matrices and hence to reformulate the
Unicity Conjecture.

\subsection{Fricke's Trace identities}\label{s:fricke} %

In this subsection we derive some of Fricke's trace identities as
needed.

\begin{prop}\label{prop:fricke}
If $X,Y \in \SLtwoC$ then %
\begin{eqnarray}
& {\rm tr}(XY)+ {\rm tr}(XY^{-1})={\rm tr}(X)\,{\rm tr}(Y); \label{eqn:XY} \\
& \hspace{-16pt} {\rm tr}^2(X)+{\rm tr}^2(Y)+{\rm tr}^2(XY)-{\rm tr}(X)\,{\rm tr}(Y)\,{\rm tr}(XY) %
=2+{\rm tr}(XYX^{-1}Y^{-1}). \label{eqn:Fv} %
\end{eqnarray}
In particular, if $X,Y \in \SLtwoC$ satisfy \,${\rm
tr}(XYX^{-1}Y^{-1})=-2$ then
\begin{eqnarray}\label{eqn:-2}
{\rm tr}^2(X)+{\rm tr}^2(Y)+{\rm tr}^2(XY)={\rm tr}(X)\,{\rm tr}(Y)\,{\rm tr}(XY). %
\end{eqnarray}
\end{prop}

\begin{pf}
These identities can be verified easily by straightforward
calculations. Here, however, we include a simpler derivation as
presented in, for instance, \cite{goldman} (see also
\cite{mumford-s-w}), which not only enables us to avoid tedious
calculations but also would led us to the
rediscovery of the identities. 

First, note that if
$Y=\Big(\footnotesize\begin{matrix}a & b \\ c & d \end{matrix}\normalsize\Big)$ then %
$Y^{-1}=\Big(\footnotesize\begin{matrix}\phantom{-}d & -b \, \\ -c & \phantom{-}a \,\end{matrix}\normalsize\Big).$ %
Hence \,${\rm tr}(Y)={\rm tr}(Y^{-1})$\, and  \,$Y+Y^{-1}={\rm
tr}(Y)\,I$, where $I$ denotes the identity matrix. Then
left-multiplying the latter equality by $X$ gives
\begin{eqnarray}\label{eqn:XtrY}
XY + XY^{-1}=X{\rm tr}(Y).
\end{eqnarray}
Taking traces on both sides of (\ref{eqn:XtrY}), we obtain identity
(\ref{eqn:XY}). As a special case, we take $X=Y$ in (\ref{eqn:XY})
to get
\begin{eqnarray}
{\rm tr}(Y^2)={\rm tr}^2(Y)-2.
\end{eqnarray}

\nn Finally, by making use of identity (\ref{eqn:XY}) repeatedly, we
can calculate ${\rm tr}(XYX^{-1}Y^{-1})$ and thus obtain
(\ref{eqn:Fv}) easily as follows:
\begin{eqnarray*}
{\rm tr}(XYX^{-1}Y^{-1}) \!\!\!&=&\!\!\! {\rm tr}(X){\rm tr}(YX^{-1}Y^{-1})-{\rm tr}(XYXY^{-1}) \\ %
\!\!\!&=&\!\!\! {\rm tr}^2(X)-\left[{\rm tr}(XY){\rm tr}(XY^{-1})-{\rm tr}(XYYX^{-1}) \right] \\ %
\!\!\!&=&\!\!\! {\rm tr}^2(X)-{\rm tr}(XY)\big[{\rm tr}(X){\rm tr}(Y)-{\rm tr}(XY)\big]+{\rm tr}(Y^2) \\ %
\!\!\!&=&\!\!\! {\rm tr}^2(X)-{\rm tr}(X){\rm tr}(Y){\rm tr}(XY)+{\rm tr}^2(XY)+{\rm tr}^2(Y)-2. %
\end{eqnarray*}

This proves Proposition \ref{prop:fricke}.
\end{pf}

\nn {\bf Remark.} Many other trace identities of Fricke for matrices
in $\SLtwoC$, though shall not be needed in this paper, have been
explored in \cite{goldman} in details.


\subsection{Markoff matrices}\label{ss:M-matrices}

Following Cohn \cite{cohn1955am} but with a different choice, we
associate a matrix in $\SLtwoZ$ to each slope $t \in \QQ$ as
follows. Initially, we set
\begin{eqnarray}\label{eqn:A=,B=}
A=\begin{pmatrix} 2 & 1 \\ 1 & 1 \end{pmatrix}, \quad %
B=\begin{pmatrix} 1 & 1 \\ 1 & 2 \end{pmatrix} %
\end{eqnarray}
and define
\begin{eqnarray}\label{eqn:M0/1=,M1/0=}
M_{\frac01}=A=\begin{pmatrix} 2 & 1 \\ 1 & 1 \end{pmatrix}, %
\quad M_{\frac10}=AB=\begin{pmatrix} 3 & 4 \\ 2 & 3 \end{pmatrix}. %
\end{eqnarray}
In general, for a Farey pair $r,s \in \hat{\mathbf Q}\cup
[0,\infty]$ with $r<s$, we set
\begin{eqnarray}\label{eqn:MrMs}
{M_{r\oplus s}=M_{r}M_{s} \, (\,\neq M_{s}M_{r}).}
\end{eqnarray}
Thus we have defined for every $t\in\hat{\mathbf Q}\cup [0,\infty]$
a {\it Markoff matrix}, {$M_t \in {\rm SL}(2,{\mathbf Z})$}, with
positive elements. As a few more examples, one finds
\begin{eqnarray*}\label{eqn:Mr=exmaples3}
{M_{\frac12}=\begin{pmatrix} 21 & 29 \\ 13 & 18 \end{pmatrix}, \quad %
M_{\frac11}=\begin{pmatrix} 8 & 11 \\ 5 & 7 \end{pmatrix}, \quad %
M_{\frac21}=\begin{pmatrix} 46 & 65 \\ 29 & 41 \end{pmatrix};}
\end{eqnarray*}
\begin{eqnarray*}\label{eqn:Mr=exmaples4}
M_{\frac13}=\begin{pmatrix} 55 & 76 \\ 34 & 47 \end{pmatrix}, \,\, %
M_{\frac23}=\begin{pmatrix} 313 & 434 \\ 194 & 269 \end{pmatrix}, \,\,%
M_{\frac32}=\begin{pmatrix} 687 & 971 \\ 433 & 612 \end{pmatrix}, \,\,%
M_{\frac31}=\begin{pmatrix} 268 & 379 \\ 169 & 239 \end{pmatrix}.
\end{eqnarray*}

\vskip 3pt

It is easy to observe that the trace of $M_t$ equals $3$ times the
$(2,1)$-element; for proof, see Proposition
\ref{prop:Mr-properties}(iii), \S \ref{s:properties}. Thus we may
write for $t\in \textstyle \QQ$
\begin{eqnarray}\label{eqn:mr}
m_t:={\rm tr}(M_t)/3.
\end{eqnarray}

Recall from \S \ref{ss:farey} that by a Farey triple $(r,t,s)$ in
$\QQ$ with $r<t<s$ we mean that $r,s\in\QQ$ are a Farey pair and
that $t=r\oplus s$.

\begin{prop} For every Farey triple $(r,t,s)$ in $\QQ$,
$(m_{r},m_{t},m_{s})$ is a Markoff triple.
\end{prop}

\begin{pf}
This follows from a simple application of identity (\ref{eqn:-2})
with $X=M_{r}$ and $Y=M_{s}$. To apply (\ref{eqn:-2}), we need to
verify that
\begin{eqnarray}\label{eqn:tr=-2}
{\rm tr}(M_r M_s M^{-1}_r M^{-1}_s)=-2
\end{eqnarray}
for every pair of Farey neighbor $r,s \in \QQ$ with $r<s$. Indeed,
since
\begin{eqnarray*}
  {\rm tr}(M_r M_t M^{-1}_r M^{-1}_t)  %
= {\rm tr}(M_r M_s M^{-1}_r M^{-1}_s)  %
= {\rm tr}(M_t M_s M^{-1}_t M^{-1}_s), %
\end{eqnarray*}
it suffices to check (\ref{eqn:tr=-2}) for the initial pair
$(r,s)=(\frac01,\frac10)$. This is true because
\begin{eqnarray*}\label{eqn:trcomm}
{\rm tr}\left(M_{\frac01}M_{\frac10}M_{\frac01}^{-1}M_{\frac10}^{-1}\right) %
=\,{\rm tr}\begin{pmatrix} -7 & 6 \\ -6 & 5 \end{pmatrix}=\,-2. %
\end{eqnarray*}
Since ${\rm tr}M_{r}=3m_{r}$ etc., we obtain from (\ref{eqn:-2})
that
\begin{eqnarray*}
(3m_r)^2+(3m_s)^2+(3m_t)^2=(3m_{r})(3m_{s})(3m_t).
\end{eqnarray*}
This shows that $(m_{r},m_{t},m_s)$ is a Markoff triple.
\end{pf}

\subsection{Matrix form of the Unicity Conjecture}

In terms of Markoff matrices defined above, we may rephrase the
Unicity Conjecture as:

\vskip 6pt

{\sc The Unicity Conjecture (Matrix form).} {\it The traces of
Markoff matrices $M_r$, $r\in \hat{\mathbf Q} \cap [0,\infty]$ are
all distinct.}


\section{Properties of Markoff matrices}\label{s:properties}

\nn The Markoff matrices defined in \S \ref{s:M-matrices} possess
certain nice properties which can be easily observed by inspecting
just a few examples.

\subsection{Elements of a single Markoff matrix} In a Markoff
matrix, we have

\begin{prop}\label{prop:Mr-properties}
For $t\in \hat{\mathbf Q} \cap [0,\infty]$, let $M_{t}=$
$\Big(\footnotesize\begin{matrix} a & b \\ c &
d\end{matrix}\normalsize\Big)$ be the Markoff matrix defined above.
Then \,{\rm(i)} $c \le d \le a \le b$; \,{\rm(ii)} $3a \ge 2b$, $3c
\ge 2d$; and \,{\rm(iii)} $a+d=3c$. Moreover, the inequalities in
{\rm(i)} and {\rm(ii)} are all strict when $t \neq 0,\infty$.
\end{prop}

\begin{pf} We prove (i)--(iii) by induction on the Farey level of
$t$. The conclusions (i)--(iii) are readily seen to be true for
$t\in \hat{\mathbf Q}\cap[0,\infty]$ of Farey level up to $1$, that
is, for $r=\frac01,\frac11,\frac10$.
Now suppose $t \in \QQ$ has Farey level at east $2$. As pointed out
in \S \ref{ss:farey}, there exists a unique Farey pair $r,s \in \QQ$
with $r<s$, such that $t=r\oplus s$. In particular, $r$ and $s$ have
smaller Farey levels. Let
\begin{eqnarray}\label{eqn:Mr=,Ms=}
{M_{r}=\begin{pmatrix} a & b \\ c & d \end{pmatrix}, \quad %
M_{s}=\begin{pmatrix} x & y \\ z & w \end{pmatrix}}.
\end{eqnarray}
Then, by definition,
\begin{eqnarray}\label{eqn:Mt=}
{M_{t}=M_{r}M_{s}=\begin{pmatrix} ax+bz & ay+bw \\
cx+dz & cy+dw \end{pmatrix}.} %
\end{eqnarray}

\nn To complete the inductive step, we proceed to prove (ii), (iii)
and (i) in this order.

\vskip 6pt

\nn {\sc Proof of (ii) for the inductive step}: \,It suffices to
observe that
\begin{eqnarray}\label{eqn:ii}
{\frac{y}{x} < \frac{ay+bw}{ax+bz} < \frac{cy+dw}{cx+dz} < \frac{w}{z} \le \frac32.} %
\end{eqnarray}

\vskip 3pt

\nn {\sc Proof of (iii) for the inductive step}: We need to show
\begin{eqnarray}\label{eqn:a+d=3c}
{(ax+bz) + (cy+dw) = 3(cx+dz).}
\end{eqnarray}
The inductive hypothesis gives
\begin{eqnarray}\label{eqn:induction iii}
{a+d=3c, \quad x+w=3z.}
\end{eqnarray}
Thus (\ref{eqn:a+d=3c}) is equivalent to
\begin{eqnarray}\label{eqn:2dx=}
{2dx=bz+cy.}
\end{eqnarray}

\nn There are two possibilities: the denominator (or numerator) of
$r$ is less or greater than that of $s$. Accordingly, we have
$s=r\oplus t'$ or $r=t'\oplus s$, where $t'\in \hat{\mathbf Q} \cap
[0,\infty]$ has Farey level lower than the maximum of those of $r$
and $s$. In the case where $s=r\oplus t'$ we have
\begin{eqnarray}\label{eqn:Mt'=}
{M_{t'}=M_{r}^{-1}M_{s}=\begin{pmatrix} d & -b \\ -c & a \end{pmatrix} %
\begin{pmatrix} x & y \\ z & w \end{pmatrix}=\begin{pmatrix} dx-bz & dy-bw \\
-cx+az & -cy+aw \end{pmatrix}.} %
\end{eqnarray}
Now the inductive hypothesis on $M_{t'}$ yields
\begin{eqnarray}\label{eqn:hypothesis on Mt'}
{(dx-bz) + (-cy+aw) = 3(-cx+az)} %
\end{eqnarray}
which is, by (\ref{eqn:induction iii}), equivalent to
(\ref{eqn:2dx=}). The proof for the other case is entirely similar.

\vskip 6pt

\nn {\sc Proof of (i) for the inductive step}: The first and the
last of the three inequalities in (i), that is,
\begin{eqnarray*}
ax+bz < ay+bw \quad\quad \text{and} \quad\quad  cx+dz < cy+dw
\end{eqnarray*}
follow easily from the inductive hypothesis $x \le y, \, z \le w$,
of which at least one inequality is strict. It remains to prove
\begin{eqnarray}\label{eqn:a>d}
{ax+bz > cy+dw.}
\end{eqnarray}
By (\ref{eqn:a+d=3c}), this is equivalent to
\begin{eqnarray*}
{3(cx+dz) > 2(cy+dw)}
\end{eqnarray*}
which is true since we have from the inductive hypothesis that $3x
\ge 2y$ and $3z \ge 2w$, and at least one of these two inequalities
is strict.

This finishes the inductive step as well as the proof of Proposition
\ref{prop:Mr-properties}.
\end{pf}

\subsection{Alternative proof of Lemma \ref{lem:um}}\label{ss:pfum} %

By Proposition \ref{prop:Mr-properties}, every Markoff matrix $M_t$,
$t \in \QQ$ is of the form
\begin{eqnarray*}
M_t=\begin{pmatrix} 2m_t-u & \ast \\ m_t & m_t+u \end{pmatrix},
\end{eqnarray*}
where $0 \le u \le m_t/2$. Now $\det M_t=1$ implies that $u^2+1
\equiv 0 \;({\rm mod}\; m_t)$. By the definition of $u_t$, we have
$u=u_t$. Thus we obtain

\begin{prop}\label{prop:Mt} %
For every $t \in \QQ$, the Markoff matrix
\begin{eqnarray}\label{eqn:Mt}
M_t=\begin{pmatrix} 2m_t-u_t & 2m_t+u_t-v_t \\ m_t &
m_t+u_t\end{pmatrix}.
\end{eqnarray}
\end{prop}

Using (\ref{eqn:Mt}), we can give an alternative proof of Lemma
\ref{lem:um}.

\vskip 6pt

\nn {\sc Alternative Proof of Lemma \ref{lem:um}.}\; We obtain from
$M_rM_s=M_t$ that
\begin{eqnarray*}
M_s=M_r^{-1}M_t=%
\begin{pmatrix} \ast & \ast \; \\ u_tm_r-u_rm_t & \ast \; \end{pmatrix}. %
\end{eqnarray*}
Equating the $(2,1)$-elements then gives the first equality in
(\ref{eqn:um}). The second equality in (\ref{eqn:um}) follows
similarly from $M_r=M_tM_s^{-1}$. \qed

\subsection{Monotonicity of the index of a Markoff matrix}

It is convenient to introduce an index
\begin{eqnarray}\label{eqn:rho}
{\varrho (M_t) := \frac{a}{c}}
\end{eqnarray}
for every Markoff matrix $M_t=$ %
$\Big(\,{\footnotesize \begin{matrix} a & b \\ c & d \end{matrix}\normalsize\,\Big)}$, %
where $r\in \QQ$.

We then have the following monotonicity of the index of a Markoff
matrix with respect to its slope. This follows readily from Lemma
\ref{lem:increasing} and (\ref{eqn:Mt}). However, we choose to
include a direct simple proof which in turn gives an alternative
proof of Lemma \ref{lem:increasing}.

\begin{prop}\label{prop:rhoM-properties} %
Suppose $t_1, t_2 \in \QQ$ where $t_1<t_2$. Then {$\varrho (M_{t_1})
> \varrho (M_{t_2})$.}
\end{prop}

\begin{pf}
We proceed to prove this proposition by induction on the maximum of
the Farey levels of $t_1$ and $t_2$.
First, the conclusion is true for $t_1, t_2$ both having Farey level
$0$, that is, $t_1=\frac01$ and $t_2=\frac10$, since
$\rho(M_{\frac01})=2/1$ and $\rho(M_{\frac10})=3/2$.

By the process of constructing all the rationals in $[0,\infty]$ by
recursively making Farey sums as described in \S \ref{ss:farey}, we
only need to prove the conclusion locally, that is, suppose it is
true for a Farey pair $r,s\in \QQ$ with $r<s$ and show
\begin{eqnarray}\label{eqn:rts}
\varrho (M_{r})>\varrho (M_{t})>\varrho (M_{s}),
\end{eqnarray} %

\nn where we have written $t:=r\oplus s$.
Suppose $M_r$, $M_s$ are given by (\ref{eqn:Mr=,Ms=}). Then $M_t$ is
given by (\ref{eqn:Mt=}). By our inductive hypothesis, $\varrho(M_r)
> \varrho(M_s)$, that is, $a/c > x/z$.
Then (\ref{eqn:rts}) is equivalent to
\begin{eqnarray}\label{eqn:a/c>rho>x/z}
{\frac{a}{c} > \frac{ax+bz}{cx+dz} > \frac{x}{z}.}
\end{eqnarray}

The first inequality in (\ref{eqn:a/c>rho>x/z}) follows easily from
the fact that $a/c>b/d$ (since $ad-bc=1$). The second is equivalent,
by Proposition \ref{prop:Mr-properties} (iii), to the inequality
\begin{eqnarray*}
{\frac{cy+dw}{cx+dz} < \frac{w}{z},}
\end{eqnarray*}
which is true since $y/x<w/z$ (as $xw-yz=1$).
This finishes the inductive step as well as the proof of Proposition
\ref{prop:rhoM-properties}.
\end{pf}

\vskip 6pt

As an immediate corollary of Proposition \ref{prop:rhoM-properties},
we obtain that

\begin{prop}\label{prop:disctinct} %
The Markoff matrices $M_t$, $t\in \QQ$ are all distinct.
\end{prop}

\nn {\bf Remark.}\, There are other choices in the definition of
Markoff matrices, such as
\begin{eqnarray}\label{eqn:new M0/1=,M1/0=}
M_{\frac01}=\begin{pmatrix} 2 & 1 \\ 1 & 1 \end{pmatrix}, %
\quad M_{\frac10}=\begin{pmatrix} 5 & 2 \\ 2 & 1 \end{pmatrix}. %
\end{eqnarray}
For this choice we have for all $t\in\QQ$
\begin{eqnarray}\label{eqn:newMt}
M_t=\begin{pmatrix} 2m_t+u_t & m_t \\ 2m_t-u_t-v_t &
m_t-u_t\end{pmatrix}.
\end{eqnarray}

\section{Further remarks}\label{s:rmk}

\nn In this section we make further remarks to give a geometric
explanation of the Markoff numbers and related numbers.

\subsection{Once-cusped hyperbolic torus}\label{ss:torus}

It is known from Cohn's work \cite{cohn1955am} that the Markoff
numbers correspond to the simple closed geodesics on a special
hyperbolic torus with a single cusp. (See \cite{series1985mi} an
exposition of the background.) Specifically, let $A, B \in \SLtwoR$
be given as in \S \ref{ss:M-matrices} and let $\langle A,B \rangle
\subset \SLtwoR$ be the subgroup generated by $A$ and $B$. Then
$\langle A,B \rangle$ is a Fuchsian group and $T:={\mathbb H}^2 /
\langle A,B \rangle$ is once-cusped hyperbolic torus, where
${\mathbb H}^2$ is the upper half-plane model of the hyperbolic
plane. Note that the axes of the M\"{o}bius transformations $A$, $B$
and $AB$ project onto simple closed curves on $T$. Assign to these
three simple closed curves on $T$ the slopes $\frac01$,
$\frac{-1}{1}$ and $\frac10$ respectively, and consider all the
simple closed cures $\gamma_t$ on $T$ of slopes $t\in\QQ$. Let the
hyperbolic length of $\gamma_t$ be $l_t$. Then we have the relation
$$ 3m_t=2\cosh (l_t/2). $$
Hence the Unicity Conjecture is actually a conjecture about the
uniqueness of lengths of certain simple closed geodesics on the
specific hyperbolic torus $T$.

\subsection{McShane identity}\label{ss:mcshane}

For a Farey triple $(r,t,s)$ in $\QQ$ with $r<t<s$, the quantity
$\displaystyle\frac{m_{t'}}{m_rm_s}$ (where $m_{t'}=3m_rm_r-m_t$)
appeared in (\ref{eqn:diff}) has nice geometric meanings. In
particular, it leads naturally to the interesting McShane identity;
see \cite{bowditch1998plms} (Theorem 3 and Proposition 3.13 there)
for details.

\subsection{Exceptional vector bundles on $\CPtwo$}

In an unexpected way, the Markoff numbers also appear as the ranks
of the exceptional vector bundles on $\CPtwo$, as explained by
Rudakov \cite{rudakov1988}. In this context, the quantity $u/m$ is
the ``slopes'' the corresponding vector bundles, with $u$ the first
Chern class (which is an integer in this case).

\end{document}